\documentclass[letterpaper, 10 pt, conference]{ieeeconf}  

\IEEEoverridecommandlockouts                   
\newcommand{\keepB}{true}
\newcommand{\revise}[1]{#1}

\usepackage[utf8]{inputenc}
\usepackage[T1]{fontenc}
\usepackage{graphicx}

\usepackage{float}
\usepackage{amsmath}
\usepackage{amssymb}

\usepackage{mathrsfs}  
\usepackage{amsthm}
\usepackage{mathtools}
\usepackage[noend,ruled,vlined,algo2e]{algorithm2e}

\usepackage{dsfont}
\usepackage{enumerate}
\usepackage[tight]{subfigure}
\usepackage{multirow}
\usepackage{bm}
\usepackage{color}

\usepackage[hidelinks]{hyperref}\usepackage{cleveref} 

\newcommand{\modif}[1]{#1}

\newcommand{\vect}[1]{#1} 
\newcommand{\m}[1]{\mathbf{#1}} 

\newcommand{\R}{\mathbb{R}}
\newcommand{\N}{\mathbb{N}}

\newcommand{\Z}{\mathbb{Z}}

\newcommand{\defEqual}{\coloneqq}
\newcommand{\argmin}{\operatorname{argmin}} 
 
\newcommand{\vol}{\operatorname{vol}} 
\newcommand{\pdc}[1]{\mathbb{S}_{+}^{#1}}

\newcommand{\co}{\operatorname{co}}
\newcommand{\myemptyset}{\varnothing}
\newcommand{\cheb}{\operatorname{CC}}
\newcommand{\FXU}{\mathscr{F}}

\newtheorem{theorem}{Theorem}
\newtheorem{lemma}{Lemma}
\newtheorem{corollary}{Corollary}
\newtheorem{proposition}{Proposition}
\newtheorem{definition}{Definition}

\newcommand{\DotSym}{\bullet} 
\newcommand{\ellipsoid}[2]{\operatorname{E}(#1,#2)} 
\newcommand{\ball}[2]{\operatorname{B}(#1,#2)} 
\newcommand{\hyperrectangle}[2]{\operatorname{H}(#1,#2)}

\newcommand{\nX}{n_x} 
\newcommand{\nU}{n_u} 
\newcommand{\nW}{n_w}

\newcommand{\set}[1]{\mathcal{#1}}
\newcommand{\Sys}{\mathcal{S}}
\newcommand{\Cont}{\mathcal{C}}
\newcommand{\seq}[1]{\boldsymbol{#1}} 

\newcommand{\keep}[2]{%
    \ifthenelse{\equal{\keepB}{true}}%
    {\modif{#1}}%
    {#2}%
}

\newcommand{\keepSecond}[2]{%
    \ifthenelse{\equal{\keepB}{true}}%
    {#1}%
    {#2}%
}

\title{
Smart abstraction based on iterative cover and non-uniform cells
}

\author{Julien Calbert \and Lucas N. Egidio \and Raphaël M. Jungers
\thanks{JC is a FRIA Research Fellow.
RJ is a FNRS honorary Research Associate. This project has received funding from the European Research Council (ERC) under \emph{European Union's Horizon 2020 research and innovation programme} under grant agreement No 864017 - L2C. The results of the numerical experiments presented in~\Cref{sec:experiments} are available in the Dionysos Julia package \url{https://github.com/dionysos-dev/Dionysos.jl/releases/tag/2024_LCSS} in the subfolder \texttt{utils/CDC2024}.}
\thanks{J.~Calbert, L.~N.~Egidio and R.~M.~Jungers are with the ICTEAM,
        UCLouvain, 4 Av. G. Lema\^{i}tre, 1348 Louvain-la-Neuve, Belgium.
{\tt\small \{julien.calbert,lucas.egidio,}
{\tt\small raphael.jungers\}@uclouvain.be}}
} 

\begin{document}

\maketitle
\thispagestyle{empty}

\begin{abstract}

We propose a multi-scale approach for computing abstractions of dynamical systems, that incorporates both \emph{local} and \emph{global} optimal control to construct a \emph{goal-specific} abstraction. 
For a local optimal control problem, we not only design the controller ensuring the transition between every two subsets (\emph{cells}) of the state space but also incorporate the volume and shape of these cells into the optimization process. This integrated approach enables the design of \emph{non-uniform} cells, effectively reducing the complexity of the abstraction.
These local optimal controllers are then combined into a digraph, which is globally optimized to obtain the entire trajectory. 
The global optimizer attempts to \emph{lazily} build the abstraction along the optimal trajectory, which is less affected by an increase in the number of dimensions. Since the optimal trajectory is generally unknown in practice, we propose a methodology based on the RRT* algorithm to determine it incrementally.
\revise{Finally, we provide a tractable implementation of this algorithm for the optimal control of $L$-smooth nonlinear dynamical systems.}
\end{abstract}

\section{Introduction}
Abstraction-based techniques have been a popular approach to safety-critical control of cyber-physical systems, enabling complex specifications and dynamics to be taken into account in the control design problem. Generally, these methods involve discretizing both the state space and the input space with uniform hyperrectangles \cite{rungger2016scots}. The curse of dimensionality, however, significantly affects the uniform discretization of the entire state space due to the exponential growth of the number of states with respect to the dimension. Additionally, to account for the quantization error between the actual state and the quantized state, it is required to over-approximate the forward image of the cells under constant discretized inputs.
One of the main drawbacks of this approach is that, in the absence of incremental stability ($\delta$-GAS)~\cite{khalil}, 
over-approximation increases the level of non-determinism in the symbolic system, which could result in an intractable or even unsolvable symbolic problem.

\revise{In~\cite{hsu2019lazy}, the authors propose to adapt the size of the abstraction gradually but uniformly
over the whole state space.}
In~\cite{calbert2021alternating}, while the authors provide a strategy to lazily (i.e., postponing heavier numerical operations) build an abstraction along the optimal trajectory, the local control design is performed combinatorially on a discretized input set, which is likely a source of non-determinism.
More recently, in~\cite{egidio2022state}, while the authors propose to design local feedback controllers between nearby cells to eliminate the non-determinism in the abstraction (i.e., the abstraction is a weighted digraph instead of a hypergraph), it relies on a predefined ellipsoidal cover of the entire state space, which suffers from the curse of dimensionality.

\begin{figure}
    \centering
    \includegraphics[width=\linewidth]{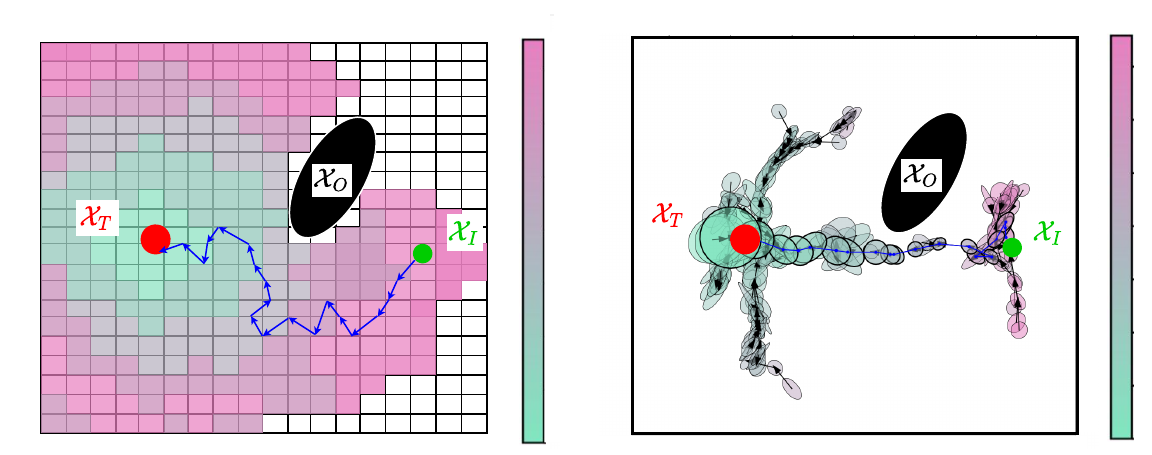}
    \caption{ 
    Comparison between the classic approach and the smart state-feedback abstraction for a planar system with state trajectory (blue line) and value function (color map) obtained for the optimal control problem of departing from $\set{X}_I$ and reaching $\set{X}_T$ while avoiding obstacles~$\set{X}_O$.
    Left: Classic abstraction approach based on predefined grid discretization of the entire state-space. Non-colored represents a region where no controller could be designed.
    Right: Smart state-feedback abstraction lazily constructed with locally optimized ellipsoids and controllers.}
    \label{fig:smart}
\end{figure}

In light of these shortcomings, we propose a new approach that leverages the incremental construction of~\cite{calbert2021alternating} and the optimal control solution proposed in~\cite{egidio2022state} by optimizing not only the transitions but also the positioning and shape of the cells along the optimal trajectory (see~\Cref{fig:smart}).
Precisely, our solution relies on a Rapidly-Exploring Random Trees (RRT) algorithm that \emph{lazily} constructs the abstraction on the basis of an \emph{ellipsoidal covering} of the state space and a finite set of \emph{local affine controllers}. 
\revise{Unlike~\cite{rungger2016scots} and~\cite{hsu2019lazy}, one of the key advantages of our approach is that, instead of discretizing the input space, we use a finite set of local affine controllers as our symbolic input set.}
\revise{ 
For $L$-smooth nonlinear systems, the combination of linearization and the Lipschitz constant allows us to locally optimize such controllers using semi-definite programming without the need for extensive discretization.
}
Moreover, the proposed approach differs from classic techniques as the partitioning is designed smartly, building the abstraction iteratively, instead of adopting a predefined uniform partition, which is suboptimal and prone to the curse of dimensionality.

The RRT algorithm~\cite{lavalle1998rapidly} is a sampling-based path planning algorithm \cite{kuffner2000rrt,lavalle2001randomized} that generates open-loop trajectories for nonlinear systems. Therefore, these techniques are not robust to perturbations and can only provide probabilistic guarantees \cite{karaman2011sampling}.  
In~\cite{manjunath2021safe}, the authors propose a safety-critical path planner based on RRT and relying on control barrier functions. 
In this same vein, the authors of~\cite{ren2019dynamic} propose an abstraction-based technique relying on the RRT algorithm limited to deterministic systems with additional stability assumptions ($\delta$-GAS). 
In~\cite{wu2022robust}, the authors propose a reachable set-based RRT planning algorithm that provides a formal guarantee in the presence of bounded perturbations, does not take into account the optimal control problem and is limited to piecewise constant controllers.
Furthermore, their framework does not allow to optimize the shape of new cells. 
Finally, the authors of~\cite{asselborn2013control} propose a forward control strategy for reachability based on ellipsoids, an affine linear approximation with an affine controller via LMIs resolution.

Our work differs from the previously mentioned works by the following points.
Firstly, our approach not only provides a controller to solve a specific control problem, but also an abstraction of the original system that can be reused as a basis for solving other optimal control problems. 
Secondly, our RRT algorithm does not act on states but on sets of states, which makes it possible to build a robust controller that provides formal guarantees even in the presence of noise. 
Thirdly, in contrast with~\cite{asselborn2013control}, our RRT algorithm operates in a backward manner. In this configuration, the target set serves as the root of the growing tree, offering the advantages of 1) maximizing the volume of newly added ellipsoids, thereby covering more of the state space with a single cell in the abstraction, and 2) effectively handling obstacles in the state space, as discussed later.

\keep{
\revise{
In~\Cref{sec:prob-statement}, we outline the fundamental concepts of our control framework. \Cref{sec:rrt} details a general method for lazily constructing an abstraction using non-uniform cells, without restricting the class of systems or cost function templates.  
From~\Cref{sec:local-transition} onward, we present a particular and optimized implementation tailored for a specific class of $L$-smooth nonlinear systems with quadratic cost functions. This approach relies on linearization to handle nonlinearities and uses the Lipschitz constant to convert the problem into a convex optimization formulation with LMIs as constraints.
}
}

\noindent\textbf{Notation:}
Given two sets $A,B$, we define a \emph{single-valued map} as $f:A\rightarrow B$, while a~\emph{set-valued map} is defined as $f:A\rightarrow 2^B$, where $2^B$ is the power set of $B$, i.e., the set of all subsets of $B$. 
The image of a subset $\Omega\subseteq A$ under $f:A\rightarrow 2^B$ is denoted $f(\Omega)$. 
\keep{Given a matrix $\m A\in \R^{m\times n}$ and a set $X\subseteq \R^n$, we define $\m A X = \{\m A \vect x \mid \vect x\in X\}.$}{}
We define $\FXU(X, U)$ as the set of \revise{single-valued functions $f:X\rightarrow U$}.
We denote by $\pdc{n}$ the set of positive definite matrices of dimension $n$. Also, $\m A\succ \m 0$  represents that  $\m A\in\pdc{n}$ and $\m A \succeq \m 0$, that $\m A\in \partial \pdc{n}$, i.e., $\m A$ is positive semidefinite. 
The function $\|\cdot\|:\R^n\rightarrow \R$ denotes the euclidean norm. 
The Minkowski sum of two sets $A,B\subset \R^n$ is $A\oplus B = \{a+b\mid a\in A,b\in B\}$.
Given a bounded set $A$, we define the \emph{Chebyshev center} of $A$ as $\cheb(A) =\argmin_{\vect x\in A} \max_{\vect y\in A}\ \|\vect x-\vect y\|$.
%
An $n$-dimensional hyperrectangle of \emph{center} $\vect c\in\R^n$ and \emph{half-lengths} $\vect h\in \R^n_{+}$ is denoted as 
$\hyperrectangle{\vect c}{\vect h} \coloneqq \{\vect x\in\R^n \mid |x_i-c_i|\le h_i \text{ for }i=1,\ldots,n\}.$
An $n$-ellipsoid with \emph{center} $\vect{c}\in\R^n$ and \emph{shape} defined by $\m P\in \pdc{n}$ is denoted as $\ellipsoid{\vect c}{\m P} \coloneqq  \{\vect x\in\R^n\mid (\vect x-\vect c)^\top \m P(\vect x-\vect c) \le 1\}$. The $n$-dimensional Euclidean ball of radius $r>0$ of center $\vect c$ is denoted as $\ball{{\vect c}}{r}\coloneqq\ellipsoid{{\vect c}}{r^{-{2}}{\m I}_n}$.

\section{Problem formulation}\label{sec:prob-statement}

To introduce the symbolic control formalism needed to support the correctness of the proposed strategy, some definitions of the control framework are required.

\begin{definition}\label{def:sys}
A \emph{transition control system} is a tuple $\Sys \coloneqq (\set{X}, \set{U}, F)$ where $\set{X}\subseteq \R^{\nX}$ and $\set{U}\subseteq \R^{\nU}$ are respectively the set of states and inputs and the set-valued map $F:\set{X}\times \set{U}\rightarrow 2^{\set{X}}$ where $F(x,u)$ give the set of states that may be reached from a given state $\vect x$ under a given input~$\vect u$.
\end{definition}
\keep{
The use of a set-valued map to describe the transition map of a system allows us to model perturbations and any kind of non-determinism in a common formalism.}{}

A tuple $(\seq x, \seq u)\in \set{X}^{[0, T[}\times \set{U}^{[0, T-1[}$ is a \emph{trajectory} of length $T\in\N\cup \{\infty\}$ starting at $x(0)\in\set{X}$ of the system $\Sys=(\set{X},\set{U}, F)$ if $\forall k\in [0,T-1[: u(k)\in \set{U}(x(k))$ and $x(k+1) \in F(x(k), u(k))$.
The set of trajectories of $\Sys$ is called the \emph{behavior} of $\Sys$, denoted $\set{B}(\Sys)$.

We consider \emph{static} controllers, where the set of control inputs enabled at a given state depends only on that state.
\begin{definition}
    We define a static \emph{controller} for a system $\Sys = (\set{X},\set{U}, F)$ valid on a subset $\set{X}_v\subseteq \set{X}$ as a \revise{single}-valued map $\Cont :\set{X}_v\rightarrow \revise{\set{U}}$ such that $\forall x\in \set{X}_v:\Cont(x)\revise{\in} \set{U}(x)$.
    %
    We define the \emph{controlled system}, denoted as $\Cont\times \Sys$, as the transition system characterized by the tuple $(\set{X}_v,\set{U}, F_{\Cont})$ where
    $x'\in F_{\Cont}(x, u) \Leftrightarrow (u\revise{=} \Cont(x)\land x'\in F(x,u))$.
\end{definition}


\keep{
\revise{In this paper, we concentrate on reach-avoid specifications as an illustrative example for the sake of clarity. However, the proposed approach can be applied to enforce a broader range of specifications.}
}{}

Given a system $\Sys=(\set{X}, \set{U}, F)$ and sets $\set{X}_I, \set{X}_T, \set{X}_O \subseteq \set{X}$, a \emph{reach-avoid} specification is defined as
\begin{align}
\Sigma = \{ &(\seq x,\seq u) \in (\set{X} \times \set{U})^\infty \mid x(0) \in \set{X}_I \Rightarrow \exists N \in \Z_+ : \nonumber \\
& \left(x(N) \in \set{X}_T \land \forall k \in [0,N[:x(k)\notin \set{X}_O\right) \}, \label{eq:reach-avoid-spec}
\end{align}
which enforces that all states in the initial set $\set{X}_I$ will reach the target $\set{X}_T$ in finite time while avoiding obstacles in~$\set{X}_O$. We use the abbreviated notation $\Sigma=[\set{X}_I,\set{X}_T,\set{X}_O]$ to denote the specification~\eqref{eq:reach-avoid-spec}.
A system $\Sys$ together with a specification $\Sigma$ constitute a \emph{control problem} $(\Sys,\Sigma)$.
Additionally, a controller $\Cont$ is said to \emph{solve} the control problem $(\Sys,\Sigma)$ if $\set{X}_I\subseteq \set{X}_v$ and $\set{B}(\Cont\times \Sys)\subseteq \Sigma$. 
We consider an \emph{optimal control problem} whose goal is to design a controller $\Cont$ enforcing reach-avoid specification $\Sigma$ while minimizing \revise{a worst-case \emph{cost function}}
\revise{
    \begin{equation}\label{eq:cost_controller}
        \mathcal{L}(\Cont) = \sup_{x_0\in \set{X}_I} \ l_{\Cont}(x_0)
    \end{equation}
    where $l_{\Cont}(x) = 0$ for $x\in\set{X}_T$ and
     \begin{align*}
         l_{\Cont}(x) &= \mathcal{J}(x, \Cont(x)) + \sup_{x'\in F(x, \Cont(x))} \ l_{\Cont}(x') \text{ for }x\in \set{X}_v\setminus \set{X}_T,
     \end{align*}
     where $\set{J}:\set{X}\times\set{U}\rightarrow\R_+$ is a given \emph{stage cost} function.
}

In the context of optimal control, we recall the notion of \emph{value function} 
\keep{(\cite[Definition 7]{legat2021abstraction})}{\revise{(\cite[Definition 3]{egidio2022state})}},
which provides a guaranteed cost for any closed-loop trajectory starting in a given subset $\set{X}_v\subseteq \set{X}$ which decreases with each time step.
\begin{definition}
    A function $v:\set{X}\rightarrow \R$ is a \emph{value function} with stage cost $\set{J}(x,u)$ for system $\Sys= (\set{X},\set{U},F)$ in $\set{X}_v\subseteq \set{X}$ if $v(x)$ is bounded from below within $\set{X}_v$ and for all $x\in\set{X}_v$ there exists $u\in\set{U}(x)$ fulfilling the Bellman inequality
    \begin{equation}
        v(x)\ge \set{J}(x,u) + \sup_{x'\in F(x,u)}\ v(x').
    \end{equation}
\end{definition}

Our abstraction-based approach relies on the notion of~\emph{state-feedback abstraction}~\cite[Definition 1]{egidio2022state} which is a specific instance of~\emph{alternating simulation relation}~\cite[Definition 4.19]{tabuada2009verification} (see~\cite[Lemma 1]{egidio2022state}).
\begin{definition}
    Consider a system $\Sys=(\set{X},\set{U}, F)$. 
    A \emph{state-feedback abstraction} $\widetilde{\Sys}$ of $\Sys$ is a system $\widetilde{\Sys} = (\widetilde{\set{X}}, \widetilde{\set{U}}, \widetilde{F})$ that satisfies the following conditions: $\widetilde{\set{X}}\subseteq 2^\set{X}$, $\widetilde{\set{U}} \subseteq \FXU(\set{X}, \set{U})$ and 
    \keepSecond{
    $$\forall \set{\xi}\in \widetilde{\set{X}}, \ \forall \kappa\in \widetilde{\set{U}}: \ \widetilde{F}(\set{\xi}, \kappa) = \{\set{\xi}_+\}$$}
    {$\forall \set{\xi}\in \widetilde{\set{X}}, \ \forall \kappa\in \widetilde{\set{U}}: \ \widetilde{F}(\set{\xi}, \kappa) = \{\set{\xi}_+\}$}
    where     \begin{equation}\label{eq:state_feedback_abstraction}
    \set{\xi}_+ \supseteq 
    \revise{\{x'\in F(x, \kappa(x)) \mid x\in \set{\xi}\}}.
    \end{equation}
\end{definition}

A value function $v:\set{X}\rightarrow \R$ with stage cost $\set{J}$ for $\Sys$ can be derived from a value function $\widetilde{v}:\widetilde{\set{X}}\rightarrow \R$ for its state-feedback abstraction $\widetilde{\Sys}$~\cite[Theorem 1]{egidio2022state} 
\begin{equation}\label{eq:lyapunov_refinement}
    v(x) = \min\{\  \widetilde{v}(\xi)~ : ~ \xi\in \widetilde{\set{X}},  \ x\in\xi\}
\end{equation}
with any cost function $\widetilde{\set{J}}$ that verifies 
\begin{equation}\label{eq:asbtract_stage_cost}
    \widetilde{\set{J}}(\xi,\kappa)\ge \set{J}(x,\revise{\kappa(x)}),\ \forall x\in \xi.
\end{equation}

 Based on the fact that any state-feedback abstraction $\widetilde{\Sys}$ is a deterministic system by definition and that the set $\widetilde{\set{U}}$ contains static (memoryless) state-feedback controllers, the state-feedback abstraction allows a specific concretization scheme, presented in the following proposition.
\begin{proposition}\label{prop:concretization_scheme}
    Let $\Sys$ be a system and $\widetilde{\Sys}$ a corresponding state-feedback abstraction. 
    Additionally, consider specifications $\Sigma=[\set{X}_I, \set{X}_T,\set{X}_O]$ with stage cost $\mathcal{J}$ for $\Sys$ and $\widetilde{\Sigma}=[\set{\xi}_I, \set{\xi}_T,\set{\xi}_O]$ with stage cost $\widetilde{\mathcal{J}}$ for $\widetilde{\Sys}$ such that 
    \begin{equation}\label{eq:abstract_spec}
        \set{X}_I\subseteq \set{\xi}_I, \ \set{\xi}_T\subseteq \set{X}_T, \ \set{X}_O\subseteq \set{\xi}_{O}
    \end{equation}
    and satisfying~\eqref{eq:asbtract_stage_cost}.
    If $\widetilde{\Cont}:\widetilde{\set{X}}_v\rightarrow \widetilde{\set{U}}$ with $\widetilde{\set{X}}_v=\widetilde{\set{X}}$ is a controller that solves the control problem $(\widetilde{\Sys},\widetilde{\Sigma})$, then for any value function $\widetilde{v}$ of $\widetilde{\Sys}$, the controller \revise{$\Cont:\set{X}_v\rightarrow \set{U}$, defined as    \begin{equation}\label{eq:concrete_controller}
        \Cont(x) = \kappa(x), \ \kappa=\widetilde{\Cont}(\xi^*), \ \xi^* = \argmin\{\widetilde{v}(\xi)~:~\xi \in \widetilde{\set{X}}, \ x\in\xi\}
    \end{equation}} 
    with $\set{X}_v = \cup_{\xi\in \widetilde{\set{X}}_v}\ \xi$, solves the control problem $(\Sys, \Sigma)$.
    In addition, the value function $v$ of $\Sys$ defined according to~\eqref{eq:lyapunov_refinement} \revise{provides an upper bound on the cost function~\eqref{eq:cost_controller}, i.e., $v(x_0)\ge \mathcal{L}(\Cont)$ for any $x_0\in\set{X}_I$}. 
\end{proposition}
\keep{
\begin{proof}
    The proof holds from the definitions of the specifications and the state-feedback abstraction.
\end{proof}
}{}

\revise{
In the next section, we present a method to solve the following problem: to lazily build a state-feedback abstraction $\widetilde{\Sys}$ of the original system, addressing an optimal control problem with reach-avoid specification $\Sigma$ and cost function $\mathcal{L}$.
}


\section{RRT*-based abstraction}\label{sec:rrt}
The nonlinear nature of the system and of the reach-avoid specification makes it extremely difficult to solve the optimal control problem directly.
Therefore, to approximately solve this problem our abstraction-based approach will split it up into several convex subproblems.
The approach involves using a RRT (Rapidly-Exploring Random Trees) algorithm that grows a state-feedback abstraction $\widetilde{\Sys}$ of $\Sys$ with an underlying \emph{tree} structure from the target set to the initial set. The nodes of this tree represent abstract states, each of which is a \emph{set} $\xi\in 2^{\set{X}}$ and transitions from child nodes to their parents are handled by \emph{local state-feedback controllers} $\kappa\in \FXU(\set{X}, \set{U})$.
While the strategy presented is theoretically applicable to any set templates for $\widetilde{\set{X}}$ and function templates for $\widetilde{\set{U}}$, we propose a practical implementation (\Cref{algo:RRTAbstraction}) using \emph{ellipsoids} and \emph{affine functions} templates, respectively, to take advantage of the power of LMIs.
The functions that implement~\Cref{algo:RRTAbstraction} are described as follows:
\begin{itemize}
    \item \texttt{getAbsSpecification($\Sigma$)}
    returns abstract specification, i.e., \keep{conservative}{} ellipsoids $\xi_I,\xi_T,\xi_O$ satisfying~\eqref{eq:abstract_spec}.
    \item \texttt{getConState($\set{X}$)} returns a candidate \keep{concrete}{} state $c\in\set{X}$.
    \item \texttt{getKClosestAbsStates($\widetilde{\set{X}}, \ \set{A}, \ K$)} returns the $K$ abstract states of $\widetilde{\set{X}}$ closest to the set $\set{A}$ according to Euclidean distance.
    \item \texttt{solveLocalProblem($c,\ \xi_+$)} returns $\set{\xi}=\ellipsoid{c}{\m P}$, a controller $\kappa$ such that 
    $\forall \vect x\in\set{\xi}: \ F(\vect x,\kappa(x))\subseteq \set{\xi}_+$ and $\kappa(\xi)\subseteq\set{U}$, and an upper bound $\widetilde{\set{J}}$ on the transition cost from $\xi$ to $\xi_+$, i.e.,
    $\tilde{\set{J}}\ge \max_{x\in\xi}\ \set{J}(x,\kappa(x)).$
    \item \texttt{handleObstacles($\xi,\ \xi_O$)} given $\xi= \ellipsoid{c}{\m P}$, returns a new ellipsoid $\set{\xi}_s = \ellipsoid{\vect c}{\gamma \m P}$ where $\gamma\ge 1$ is the smallest value such that $\set{\xi}_s$ does not intersect any obstacles, i.e. $\set{\xi}_s\cap \set{\xi}_{O} = \myemptyset$. 
    \item \texttt{getAbsValueFunction($\widetilde{\Sys},\ \xi_T$)} returns a value function $\widetilde{v}$ for $\widetilde{\Sys}$ such that $\widetilde{v}(\xi_T) = 0$ and which is strictly positive for the other abstract states.
    \keep{
    \item \texttt{newTransition($\xi',\ \xi$)} returns a controller $\kappa$ that maps points from a given ellipsoid $\set{\xi}'$ to another given ellipsoid $\set{\xi}$, and the associated transition costs $\widetilde{\set{J}}'$.}{} 
\end{itemize}

\SetNlSty{textbf}{}{ }\IncMargin{2em}
\begin{algorithm2e}
\small{
  \SetAlgoLined
     \nl $\widetilde{\Sigma}$ $\gets$ \footnotesize{\texttt{getAbsSpecification}} \small ($\Sigma$) \;
    \nl $\set{\xi}_I,\set{\xi}_T, \set{\xi}_O$ $\gets$ $\widetilde{\Sigma}$\;
   \nl $\widetilde{\set{X}}$ $\gets$ $\{\set{\xi}_T\}$ \;
   \nl \label{line:while} \While{\normalfont not exists $\xi\in\widetilde{\set{X}}:\ \xi_I\subseteq \xi$}{
       \nl  $\vect c$ $\gets$  \texttt{getConState($\set{X}$)} \label{line:getConState}\; 
       \nl $\set{\xi}_+$ $\gets$ \footnotesize{\texttt{getKClosestAbsStates}}\small($\widetilde{\set{X}}$, $\{\vect c\}$, $K$$=$$1$)\;
       \nl \label{line:local_problem} $\set{\xi}$, $\kappa$, $\widetilde{\set{J}}$ $\gets$ \texttt{solveLocalProblem($\vect c$, $\set{\xi}_+$)}\;
       \nl \textbf{if} not \textit{feasible} \textbf{then} \Return to~\Cref{line:getConState}. \;
       \nl $\set{\xi}_s$ $\gets$ \texttt{handleObstacles($\set{\xi}$, $\set{\xi}_O$)}\;
       \nl $\widetilde{\set{X}}, \ \widetilde{\set{U}} \ $ $\gets$ 
       $\widetilde{\set{X}}\cup \{\xi_s\},\ \widetilde{\set{U}}\cup \{\kappa\}$\;
       \keep{
       \nl $\widetilde{F}(\xi_s, \kappa) = \{\xi_+\}$\;
       \nl $\widetilde{\set{J}}(\xi_s, \kappa) = \widetilde{\set{J}}$\;}{
        \nl $\widetilde{F}(\xi_s, \kappa),\ \widetilde{\set{J}}(\xi_s, \kappa) \gets \{\xi_+\},\ \widetilde{\set{J}}$\;
       }
       \nl $\widetilde{\Sys}$ $\gets$ $(\widetilde{\set{X}}, \widetilde{\set{U}}, \widetilde{F})$\;
       \nl $\widetilde{v}$ $\gets$ \texttt{getAbsValueFunction}($\widetilde{\Sys}$, $\xi_T$)\;
       \keep{\nl \texttt{improveAbs($\widetilde{\Sys}$, $\widetilde{v}$,  $\set{\xi}_s$, $\widetilde{\set{J}}$)}\;}{}
   
   }
   \nl\Return{\normalfont $\widetilde{\Sys},\ \widetilde{v}$\;}
  \caption{
    Lazy construction of an ellipsoid-based abstraction $\widetilde{\Sys}$ of $\Sys=(\set{X},\set{U},F)$ and of an abstract value function $\widetilde{v}$ that solves the specification $\Sigma = [\set{X}_I, \set{X}_T, \set{X}_O]$. 
  }
  \label{algo:RRTAbstraction} }
\end{algorithm2e}\DecMargin{2em}

\keep{
\SetNlSty{textbf}{}{ }\IncMargin{2em}
\begin{algorithm2e}
\small{
  \SetAlgoLined
  \SetKwFunction{FMain}{improveAbs}
  \SetKwProg{Fn}{function}{:}{}
  \Fn{\FMain{$\widetilde{\Sys}$, $\widetilde{v}$,  $\set{\xi}$, $\widetilde{\set{J}}$}}{
   \nl $\set{\xi}_{\text{list}}$ $\gets$  \footnotesize{\texttt{getKClosestAbsStates}}\small ($\widetilde{\set{X}}$, $\set{\xi}$, $K$$\ge$$1$)\;
   \nl \For{\normalfont $\set{\xi}'\in \set{\xi}_{\text{list}}$}{
        \nl $\kappa'$, $\widetilde{\set{J}}'$ $\gets$  \texttt{newTransition($\set{\xi}'$, $\set{\xi}$)}\;
        \nl \uIf{\normalfont $\widetilde{\set{J}}' + \widetilde{v}(\xi) < \widetilde{v}(\xi')$}
        {
            \nl $\widetilde{\set{U}}$ $\gets$ 
            $\widetilde{\set{U}}\cup \{\kappa'\}$\;
            \nl $\widetilde{F}(\xi', \kappa') = \{\xi\}$\;
            \nl $\widetilde{\set{J}}(\xi', \kappa') = \widetilde{\set{J}}'$\; 
        }
        \textbf{end}
    }
    }
     \textbf{end}
  \caption{
    Update of the abstraction for the RRT* variant.%
  }
  \label{algo:RRT*function} }
\end{algorithm2e}\DecMargin{2em}
}
{}

As mentioned in the introduction, our RRT algorithm works backwards in the sense that the abstraction is initialized with the target set $\set{\xi}_T$, and when designing a new transition, the optimized cells are preceding ones (in terms of the controlled system trajectory).
This allows us to maximize the volume of the preceding ellipsoid, covering a larger portion of the state space with a single abstract state—a feature not achievable in the forward approach, which focuses on minimizing the volume of the destination ellipsoid~\cite{asselborn2013control}.
In addition, the backward approach makes it possible to decouple the design of the preceding ellipsoid $\xi$ and the controller $\kappa$ from obstacle handling. 
Indeed, when a preceding ellipsoid $\xi$ returned by \texttt{solveLocalProblem} intersects an obstacle, it can be shrunk into $\xi_s\subseteq \xi$ and the same controller $\kappa$ still ensures a transition to $\xi_+$.
This decoupling is not feasible in the forward approach.

\keep{
The function \texttt{improveAbs} implements the RRT* variant. When a new abstract state $\xi_s$ is added, new transitions are computed from the closest existing abstract states to $\xi_s$ in order to eventually reduce the path cost to the target set $\xi_T$.}{}

\revise{
The function \texttt{getConState} is a \emph{heuristic} used to generate a concrete state from which a new ellipsoid is added to the abstraction. Although the specific implementation of this function does not affect the correctness of~\Cref{algo:RRTAbstraction}, it can be used to guide the exploration of the state space and the construction of the abstraction in order to reduce the number of abstract states. For example, a simple but effective strategy involves sampling points predominantly in the direction of the initial set $\set{X}_I$.
}


The following theorem guarantees the validity of the approach.
\begin{theorem}\label{th:RRTAbstraction}
Consider the system $\widetilde{\Sys} = (\widetilde{\set{X}}, \widetilde{\set{U}}, \widetilde{F})$ and its value function $\widetilde{v}$ generated by \Cref{algo:RRTAbstraction}.
Then 1) $\widetilde{\Sys}$ is a state-feedback abstraction of $\Sys$; 2) there exists a controller $\widetilde{\Cont}$ that solves $(\widetilde{\Sys},\widetilde{\Sigma})$; 3) the controller $\Cont$ defined according to~\eqref{eq:concrete_controller} solves $(\Sys, \Sigma)$; 4) the function $v$ derived from $\widetilde{v}$ according to~\eqref{eq:lyapunov_refinement} is a value function for $\Sys$ with $\set{X}_v= \cup_{\xi\in \widetilde{\set{X}}}\ \xi$.
\end{theorem}
\keep{
\begin{proof}
1) The transition map $\widetilde{F}$ is constructed according to $\widetilde{F}(\set{\xi}_s,\kappa) = \{\set{\xi}_+\}$, where \texttt{solveLocalProblem} guarantees that $\forall x\in \set{\xi}:\ F(x,\kappa(x))\subseteq \set{\xi}_+$, and \texttt{handleObstacles} ensures that $\xi_s\subseteq \xi$, which (always) implies $\forall x\in\set{\xi}_s: \ F(x,\kappa(x))\subseteq \set{\xi}_+$. Consequently, $\widetilde{F}$ satisfies the condition~\eqref{eq:state_feedback_abstraction}.
\\
\indent 2) 
Firstly, for all $\set{\xi}\in \widetilde{\set{X}}$, there exists a path from $\set{\xi}$ to $\set{\xi}_T$ since $\widetilde{\set{X}}$ is initialized with $\xi_T$ and newly added transitions are always directed to an existing abstract state in $\widetilde{\set{X}}$. In addition, these paths avoid obstacle $\xi_O$ thanks to \texttt{handleObstacles}.
Secondly, the ending condition of~\Cref{algo:RRTAbstraction} guarantees that $\xi_I\subseteq \xi$ for some $\xi\in\widetilde{\set{X}}$.
As a result, there exists a controller $\widetilde{\Cont}$ that solves $(\widetilde{\Sys}, \widetilde{\Sigma})$ where $\widetilde{\Sigma} = [\xi_I,\xi_T,\xi_O]$.
\\
\indent 3) Since $\widetilde{\Sys}$ is a state-feedback abstraction of $\Sys$ and $\widetilde{\Sigma}$ satisfies~\eqref{eq:abstract_spec} according to \texttt{getAbsSpecification}, then, by~\Cref{prop:concretization_scheme}, the controller~\eqref{eq:concrete_controller} solves $(\Sys, \Sigma)$.
\\
 \indent 4) This follows directly from the fact that the abstract stage cost function $\tilde{\set{J}}$ satisfies~\eqref{eq:asbtract_stage_cost} since it is defined by the upper bounds costs of the local transition according to \texttt{solveLocalProblem}.
\end{proof}
}{
\begin{proof}
Due to space constraints, we will demonstrate the first item and refer to the extended version of the paper for the proof of the remaining points.
1) The transition map $\widetilde{F}$ is constructed according to $\widetilde{F}(\set{\xi}_s,\kappa) = \{\set{\xi}_+\}$, where \texttt{solveLocalProblem} guarantees that $\forall x\in \set{\xi}:\ F(x,\kappa(x))\subseteq \set{\xi}_+$, and \texttt{handleObstacles} ensures that $\xi_s\subseteq \xi$, which (always) implies $\forall x\in\set{\xi}_s: \ F(x,\kappa(x))\subseteq \set{\xi}_+$. Consequently, $\widetilde{F}$ satisfies the condition~\eqref{eq:state_feedback_abstraction}.
\end{proof}
}

The abstract value function $\widetilde{v}$ is derived in \texttt{getAbsValueFunction} by computing the shortest path to $\xi_T$, which can be done efficiently since $\xi_T$ is the root of a tree.
The optimization problem in \texttt{handleObstacles} and the inclusion test in~\Cref{line:while} of~\Cref{algo:RRTAbstraction} can be efficiently addressed by solving a convex scalar optimization problem, as outlined in~\cite[Corollary 2]{calbert2023efficient} and~\cite[Algorithm 1]{calbert2023efficient}, respectively.
Finally, the function\keep{s}{} \texttt{solveLocalProblem} \keep{and \texttt{newTransition}}{} can be implemented efficiently by solving a convex optimization problem as described in the following section.

\section{Local controller design}\label{sec:local-transition}
In this section, we provide a practical and efficient implementation of the function \texttt{solveLocalProblem} in~\Cref{line:local_problem} of~\Cref{algo:RRTAbstraction} for the following general class of systems.

We consider nonlinear discrete-time system $\Sys = (\set{X}, \set{U}, F)$ with bounded disturbances, i.e.,
\begin{equation}\label{eq:bounded_disturbances}
    F(x,u) = \{f(x,u,w) \mid w\in \set{W}\}
\end{equation}
where $\set{X}\subseteq \R^{\nX}$, $\set{U}\subseteq\R^{\nU}$,
$f:\set{X}\times \set{U}\times \set{W}\rightarrow \set{X}$ is a continuous nonlinear function whose Jacobian is globally Lipschitz continuous
\keep{\footnote{Note that similar results can be adapted for the locally Lipschitz continuous case.}}{}, and $\set{W}\subseteq \R^{\nW}$ is a bounded polytopic set.
The control input takes values in the intersection of (possibly degenerated) ellipsoids
\begin{equation}\label{eq:inputSet}
\set{U} = \bigcap\limits_{k=1,\ldots,N_u} \set{U}_k
\end{equation}
with $\set{U}_k = \{\vect u\in \R^{\nU}: \|\m U_k \vect u\|\le 1\}$ for some matrix $\m U_k$ of appropriate dimensions.
The exogenous input takes values inside the convex hull of a set of points
\begin{equation}\label{eq:noiseSet}
\set{W} = \co\{\vect w_1,\ldots,\vect w_{N_w}\}
\end{equation}
with $\vect 0\in \set{W}$ and introduces a level of non-determinism to the system, which can either capture non-modeled behaviors or adversarial disturbances. 
We consider a general quadratic stage cost function of the form 
\begin{equation}\label{eq:cost_function}
    \set{J}(\vect x,\vect u) = 
    \begin{pmatrix}
    \vect x^\top & \vect u^\top & 1 
    \end{pmatrix} \m Q \begin{pmatrix}
    \vect x^\top & \vect u^\top & 1 
    \end{pmatrix}^\top
\end{equation}
defined for some given matrix $\m Q\succ\m 0$.

Therefore, the objective of this section is to design a controller $\kappa$ to map all the states of $\set{\xi}\in 2^{\set{X}}$ into $\set{\xi}_+\in 2^\set{X}$ despite exogenous noise, i.e.,
\begin{equation}\label{eq:LMI_transition_desired}
    \forall \vect x\in \set{\xi}\ \forall\vect w\in \set{W}:\ f(\vect x,\kappa(\vect x), \vect w)\in \set{\xi}_+
\end{equation}
and such that the inputs of the closed loop system lie within the admissible input set $\set{U}$, i.e.,
\begin{equation}\label{eq:LMI_input_feasibility}
\kappa(\set{\xi})\subseteq \set{U}.
\end{equation}

\subsection{Specific approximation scheme}

Since optimizing directly on the nonlinear dynamics $f$ is a challenging problem, we rely on the linearized function $\tilde{f}$ around a point $\vect{\bar{p}} = (\vect{\bar{x}}, \vect{\bar{u}}, \vect{\bar{w}})\in \set{X}\times \set{U} \times \set{W}$:
\begin{equation}\label{eq:linearization}
    \tilde{f}(\vect x,\vect u,\vect w) = \m A \vect x + \m B \vect u + \m E \vect w + \vect g
\end{equation}
\keep{with $\m A = J_{f,x}(\vect{\bar{p}}), \ \m B = J_{f,u}(\vect{\bar{p}}), \ \m E = J_{f,w}(\vect{\bar{p}}), \ \vect g = f(\vect{\bar{x}},\vect{\bar{u}},\vect{\bar{w}}) -\m A\vect{\bar{x}} - \m B\vect{\bar{u}}-\m E \vect{\bar{w}}$, where $J_{f,x}(\bar{p})$ is the Jacobian matrix of $f$ with respect to the variables $x$ evaluated at $\bar{p}$.}{
where $\vect g = f(\vect{\bar{x}},\vect{\bar{u}},\vect{\bar{w}}) -\m A\vect{\bar{x}} - \m B\vect{\bar{u}}-\m E \vect{\bar{w}}$.
}

We can derive a component-wise bound on the linearization error~\cite[Lemma~1.2.3]{nesterov2018lectures}, i.e.,
$\forall (x,u,w)\in \set{X}\times \set{U}\times \set{W}:$
\begin{equation}\label{eq:linearized_model}
    f(x,u,w) \in \{\tilde{f}(x,u,w)\} \oplus \Omega(x,u,w)
\end{equation}
where $\Omega(x,u,w) = \hyperrectangle{0}{\tfrac{1}{2}Lr_{\bar{p}}(x,u,w)^2}$ with $L\in \R^{\nX}$ representing the vector of Lipschitz constants of component functions $f_i$ and 
\begin{equation}\label{eq:distance}
r_{\vect{\bar{p}}}(\vect x,\vect u,\vect w)^2 = \|\vect x-\vect{\bar{x}}\|^2 + \|\vect u-\vect{\bar{u}}\|^2 + \|\vect w-\vect{\bar{w}}\|^2.
\end{equation}

Therefore, to enforce the condition~\eqref{eq:LMI_transition_desired}, we can impose the following stronger condition 
{\small\begin{equation}\label{eq:inclusion_linearized}
    \forall x\in \xi \ \forall w\in\set{W}:\{\tilde{f}(x,\kappa(x),w)\} \oplus \Omega(x,\kappa(x),w)\subseteq \set{\xi}_+.
\end{equation}}

In order to keep the future optimization problem tractable (convex), we use a common error bound corresponding to the "worst case" to impose the more conservative condition
\begin{equation}\label{eq:inclusion_linearized_conservative}
\forall x\in \xi \ \forall w\in\set{W}:\{\tilde{f}(x,\kappa(x),w)\} \oplus \set{H}_r\subseteq \set{\xi}_+,
\end{equation}
with $\set{H}_r\defEqual \hyperrectangle{0}{\tfrac{1}{2}L r^2}$ and 
$r \defEqual \max_{(x,w)\in \xi\times \set{W}}\ r_{\bar{p}}(x,\kappa(x),w)$.

Condition~\eqref{eq:inclusion_linearized_conservative} can be reformulated as
\begin{equation}\label{eq:minkowski-sum}
q(\set{\xi}) \oplus \m E\set{W} \oplus \set{H}_r \subseteq \set{\xi}_+
\end{equation}
where $q(\vect x) \defEqual \tilde{f}(\vect x, \kappa(\vect x), \vect{\bar{w}})$.
This provides a geometric interpretation of the different sources of conservatism. The term $q(\xi)$ corresponds to the closed-loop image of $\xi$ in the linearized model with nominal noise ($\bar{w}$), the second term $\m E\set{W}$ corresponds to the exogenous noise, and the third term $\set{H}_r$ is responsible for the linearization error.



We first propose the following lemma, which provides sufficient conditions to design $\kappa$ satisfying~\eqref{eq:LMI_transition_desired}.
\begin{lemma}\label{lem:conservatism}
Let a system~\eqref{eq:bounded_disturbances} with exogenous input set~\eqref{eq:noiseSet}, a controller $\kappa$, two sets  $\set{\xi}$ and $\set{\xi}_+$ and a linearization point~$\bar{p}$.
If the following conditions are satisfied
\begin{align}
        \set{\xi}&\subseteq \ball{\bar{x}}{\sqrt{\delta_X}}, \label{eq:LMI_state_bound}\\
        \kappa(\set{\xi})&\subseteq \ball{\vect{\bar{u}}}{\sqrt{\delta_U}}, \label{eq:LMI_input_bound}\\
        \forall \vect x\in \set{\xi}\ \forall \vect w\in\set{W}: &\ \{\tilde{f}(x,\kappa(w),w)\}\oplus \set{H}_r \subseteq \set{\xi}_+,
        \label{eq:LMI_transition}
\end{align}
with $r^2 = \delta_X + \delta_U + \delta_W$ for some $\delta_X,\delta_U\ge 0$ and $\delta_W \coloneqq \max_{\vect w\in\set{W}}\|\vect w-\vect{\bar{w}}\|^2$, then~\eqref{eq:LMI_transition_desired} holds.
\end{lemma}
\keep{
\begin{proof}
\modif{The result follows directly from this sequence of implications
\begin{equation*}\label{eq:conservatism_reformulated}
\eqref{eq:LMI_state_bound},\eqref{eq:LMI_input_bound},\eqref{eq:LMI_transition}\Rightarrow \eqref{eq:inclusion_linearized_conservative} 
 \Rightarrow \eqref{eq:inclusion_linearized}\Rightarrow \eqref{eq:LMI_transition_desired}.
\end{equation*}
}
\end{proof}
}{}

\keep{
However, the converse result does not generally hold. The conservatism arises from two sources. 
Firstly, we use the continuity property yielding the upper-bound~\eqref{eq:linearized_model}. Secondly, we consider the worst linearization error for all points of $\xi$, which corresponds to the error of the farthest point from the linearization point as shown in~\eqref{eq:LMI_state_bound} and \eqref{eq:LMI_input_bound}.}{}

The linearization point that minimizes the linearization error ($r$) is determined by $\bar{p}^* = \cheb(\set{\xi}\times \kappa(\set{\xi})\times \set{W})$, the Chebyshev center of $\set{\xi}\times \kappa(\set{\xi})\times \set{W}$.
However, given that the controller $\kappa$ is part of the optimization process, we opt for the choice:
$\bar{p} = \cheb(\set{\xi}\times \set{U}\times \set{W})$
which represents the optimal choice for minimizing the linearization error when $\kappa$ is arbitrary.


\subsection{Discretization and controller templates}
\keep{
An important classic result recalled in this section is the so-called S-procedure~\cite[Theorem~2.2]{polik2007survey}.
\begin{lemma}[S-procedure \cite{polik2007survey}]\label{lemma:Sprocedure}
Let $q_i:\R^n\rightarrow \R:\vect x \rightarrow q_i(\vect x) = \vect x^\top \m P_i\vect x + 2\vect g_i^\top\vect x + s_i$ for $i=0,1$ be two quadratic
functions and suppose that there is a $\vect{\bar{x}}$ such that $q_1(\vect{\bar{x}})<0$.
Then the following two statements are equivalent
\begin{align*}
    (i) &\ \forall \vect x\in\R^n:\ q_1(\vect x)\le 0\Rightarrow q_0(\vect x)\le 0\\
    (ii)&\ \exists \lambda\ge 0: \ \lambda \begin{pmatrix}
    \m P_1 & \vect g_1\\
    \vect g_1^\top & s_1
    \end{pmatrix}
    \succeq 
    \begin{pmatrix}
        \m P_0 & \vect g_0\\
        \vect g_0^\top & s_0
    \end{pmatrix}.
\end{align*}
\end{lemma}
}{}

In this paper, we discretize the state space using ellipsoids, i.e.,
\keepSecond{$$\set{\xi}=\ellipsoid{\vect c}{\m P}, \ \set{\xi}_+ =\ellipsoid{\vect c_+}{\m P_+}$$}{$\set{\xi}=\ellipsoid{\vect c}{\m P}, \ \set{\xi}_+ =\ellipsoid{\vect c_+}{\m P_+}$} with $c,c_+\in\R^{\nX}$ and $\m P,\m P_+\in \pdc{n}$, and we consider affine controllers of the form 
\begin{equation}\label{def:affine_controller}
    \kappa(\vect x) = \m K (\vect x-\vect c)+ \vect l
\end{equation}
where $\m K\in \R^{\nU \times \nX}$ and $\vect l\in \R^{\nU}$.

Therefore, in this setting, we obtain $\vect{\bar{x}} = \cheb(\set{\xi}) = \vect{c}$, and, for the sake of clarity throughout the rest of the paper, we assume that $\bar{w} = \cheb(\set{W}) = 0$.

The following theorem provides 
conditions based on \Cref{lem:conservatism} to guarantee the existence of a valid controller satisfying~\eqref{eq:LMI_transition_desired} and~\eqref{eq:LMI_input_feasibility}.
\begin{theorem}\label{th:feasibility_linearization}
The system \eqref{eq:bounded_disturbances} with control input set \eqref{eq:inputSet} and exogenous input set \eqref{eq:noiseSet} under the constrained affine control law \eqref{def:affine_controller} satisfies the condition that $\vect x(k+1)\in \set{\xi}_+=\ellipsoid{\vect c_+}{\m P_+}$ for all $\vect x(k)\in \set{\xi}=\ellipsoid{\vect c}{\m P}$, if, given the linearized system \eqref{eq:linearization} around the point $\vect{\bar{p}}=(\vect c,\vect{\bar{u}},\vect{0})$ with $\vect{\bar{u}}\in\set{U}$, 
there exist $\m L\in \pdc{\nX},\m F\in \R^{\nX\times \nU}$ and scalars $\delta_X\ge 0,\delta_U\ge 0,\phi\ge 0,\beta_{ij}\ge 0,\tau_k\ge 0$ such that 
\begin{align}
    & \begin{pmatrix}
        \m I & \m L\\
        \DotSym & \delta_X \m I
    \end{pmatrix}\succeq \m 0 
    \label{eq:th_LMI_state_bound}
    \\
    & \begin{pmatrix}
    \phi \m I & \vect 0 & \m F^\top \\
    \DotSym & \delta_U-\phi & (\vect l-\vect{\bar{u}})^\top \\
    \DotSym & \DotSym & \m I
    \end{pmatrix}\succeq \m 0 
    \label{eq:th_LMI_input_bound}
    \\
    & 
     \begin{pmatrix}
    \beta_{ij} \m I & \vect 0 & (\m A\m L + \m B \m F)^\top \\
    \DotSym & 1-\beta_{ij} & \vect \mu^\top + \vect V_i^\top +(\m E\vect w_j)^\top \\
    \DotSym & \DotSym & \m P_+^{-1}
    \end{pmatrix} \succeq \m 0, \nonumber\\ 
    &\text{ for } i\in \{1,\ldots,2^{\nX}\}, \ j\in \{1,\ldots,N_w\} 
    \label{eq:th_LMI_transition}
    \\
    &  
    \begin{pmatrix}
    \tau_k \m I & \vect 0 & \m F^\top \m U_k^\top \\
    \DotSym & 1-\tau_k & \vect l^\top \m U_k^\top  \\
    \DotSym & \DotSym & \m I
    \end{pmatrix}\succeq \m 0, \ \text{ for } k\in \{1,\ldots,N_u\} \label{eq:th_LMI_input_feasibility}
    \end{align}    
where $\vect \mu = \vect g+\m A\vect c+\m B\vect l-\vect c_+$, $\vect V_i$ for $i=1,\ldots,2^{\nX}$ are the vertices of $\hyperrectangle{\vect 0}{\tfrac{1}{2}\vect Lr^2}$ with $r^2 = \delta_X+\delta_U+\delta_W$
and $\delta_W = \max_{i=1,\ldots,N_w}\{\|\vect w_i\|^2\}$.
\keep{
From which we have
$$\m P = \m L^{-2},\ \m K = \m F\m L^{-1}.$$
}{
From which we have
$$\m P = \m L^{-2},\ \m K = \m F\m L^{-1}.$$
}

\end{theorem}

\keep{
\begin{proof}
Given \Cref{lem:conservatism}, it is sufficient to establish the following equivalences
    \begin{align*}
         \eqref{eq:LMI_state_bound} \Leftrightarrow \eqref{eq:th_LMI_state_bound},\ \eqref{eq:LMI_input_bound} \Leftrightarrow \eqref{eq:th_LMI_input_bound},\ \eqref{eq:LMI_transition} \Leftrightarrow \eqref{eq:th_LMI_transition},\ \eqref{eq:LMI_input_feasibility} \Leftrightarrow \eqref{eq:th_LMI_input_feasibility}.
    \end{align*}
     The closed-loop dynamics can be expressed as
        \begin{align*}
            \tilde{f}(\vect x,\kappa(\vect x),\vect w) &= \m{\bar{A}} \vect x + \vect{\bar{b}} +\m E\vect w
        \end{align*}
        with 
        $$\m{\bar{A}} = \m A +\m B\m K, \  \vect{\bar{b}} = \vect g + \m B\vect l-\m B\m K\vect c.$$
    We define the matrix $\m \Theta$ as
    $$ \m \Theta = 
    \begin{pmatrix}
        \m L & \vect c & \m 0\\
        \DotSym & 1 & \vect 0^\top\\
        \DotSym & \DotSym & \m I
    \end{pmatrix}.$$
    \\
    \indent
    Using respectively the S-procedure (\Cref{lemma:Sprocedure}), the fact that $\m P = \m L^{-2}$ with $\m L \succ \m 0$ and the Schur Complement Lemma \cite{boyd1994linear}, we have
    $$\eqref{eq:LMI_state_bound} \Leftrightarrow \delta_X \m P - \m I \succeq \m 0\Leftrightarrow \m I - \delta_X^{-1} \m L^2 \succeq \m 0\Leftrightarrow \eqref{eq:th_LMI_state_bound}.$$
    \\
    \indent
    By using the S-procedure, equation \eqref{eq:LMI_input_bound} can be expressed as the existence of $\phi\ge 0$ that satisfies the inequality
    \begin{multline*}\phi 
        \begin{pmatrix}
            \m P & -\m P\vect c\\
             \DotSym & \vect c^\top \m P\vect c-1
        \end{pmatrix} \succeq \\
        \begin{pmatrix}
            \m K^\top \m K & \m K^\top (\vect l-\vect{\bar{u}}-\m K\vect c)\\
            \DotSym & (\vect l-\vect{\bar{u}}-\m K\vect c)^\top  (\vect l-\vect{\bar{u}}-\m K\vect c) - \delta_U
        \end{pmatrix}.
     \end{multline*}
       By further algebraic manipulations and using the Schur Complement Lemma, this inequality can be written as
        $$\begin{pmatrix}
            \phi \m P & -\phi\m P\vect c & \m K^\top\\
            \DotSym & \phi(\vect c^\top \m P\vect c -1) + \delta_U & (\vect l-\vect{\bar{u}} -\m K\vect c)^\top \\
            \DotSym & \DotSym & \m I
         \end{pmatrix}\succeq \m 0.$$
        Then, by applying the congruent transformation of matrix $\m \Theta$, i.e., multiplying the above inequality to the right by $\m \Theta$ and to the left by $\m \Theta^\top$, we obtain equation \eqref{eq:th_LMI_input_bound}.
        \\
        \indent
        Since a polytope is contained in a convex set if and only if its vertices are contained in this set, equation~\eqref{eq:LMI_transition} is satisfied $\forall \vect w\in\set{W}$ if and only if it is satisfied for the vertices of $\set{W}$. Hence, \eqref{eq:LMI_transition} is equivalent~to 
        \begin{equation*}
        \forall j\in\{1,\ldots,N_w\} \ \forall \vect x \in \set{\xi}: \{\tilde{f}(\vect x,\kappa(\vect x),\vect w_j)\} \oplus \set{H}_r\subseteq \set{\xi}_+
        \end{equation*}
        where $\vect w_j$ are the vertices of $\set{W}$.
        By using the S-procedure, equation \eqref{eq:LMI_transition} can be expressed as the existence of $\beta_{ij}\ge 0$ that satisfies the inequality
        \begin{multline*}
        \small{\beta_{ij} \begin{pmatrix}
        \m P & -\m P\vect c\\
        \DotSym & \vect c^\top \m P\vect c-1
        \end{pmatrix} \succeq} \\
        \footnotesize{\begin{pmatrix}
        \m{\bar{A}}^\top\m P_+ \m{\bar{A}} & \m{\bar{A}}^\top \m P_+(\vect{\bar{b}}+\m E\vect w_j + \vect V_i -\vect c_+)\\
        \DotSym & (\vect{\bar{b}}+\m E\vect w_j + \vect V_i-\vect c_+)^\top \m P_+ (\vect{\bar{b}}+\m E\vect w_j + \vect V_i-\vect c_+) - 1
        \end{pmatrix}}
        \end{multline*}
        for $i\in\{1,\ldots,2^{\nX}\}$ and $j\in\{1,\ldots,N_w\}$.
        By further algebraic manipulations and using the Schur Complement Lemma, this inequality can be written as
        $$ 
         \small{\begin{pmatrix}
            \beta_{ij} \m P & -\beta_{ij}\m P\vect c & \m{\bar{A}}^\top\\
            \DotSym & \beta_{ij}(\vect c^\top \m P\vect c -1) + 1 & (\vect{\bar{b}}+\m E\vect w_j + \vect V_i-\vect c_+)^\top \\
            \DotSym & \DotSym & \m P_+^{-1}
         \end{pmatrix}\succeq \m 0.}$$
        Then, by applying the congruent transformation of matrix $\m \Theta$, we obtain equation \eqref{eq:th_LMI_transition}.
        \\
        \indent
        By using the S-procedure, equation \eqref{eq:LMI_input_feasibility} can be expressed as the existence of $\tau_k\ge 0$ that satisfies the inequality
        \begin{multline*}
        \tau_k \begin{pmatrix}
            \m P & -\m P\vect c\\
            \DotSym &  \vect c^\top\m P\vect c-1
        \end{pmatrix}\succeq \\
        \begin{pmatrix}
            \m U_k^\top\m K^\top \m K \m U_k & \m K^\top \m U_k^\top \m U_k(\vect l-\m K\vect c)\\
            \DotSym & (\vect l-\m K\vect c)^\top \m U_k^\top \m U_k  (l-\m K\vect c) - 1
        \end{pmatrix}
        \end{multline*}
        for $k\in\{1,\ldots,N_u\}$.
        By further algebraic manipulations and using the Schur Complement Lemma, this inequality can be written as
        $$ 
         \begin{pmatrix}
            \tau_k \m P & -\tau_k \m P\vect c & \m K^\top \m U_k^\top\\
            \DotSym & \tau_k (\vect c^\top\m P\vect c-1)+1 &  (\vect l-\m K\vect c)^\top \m U_k^\top\\
            \DotSym & \DotSym & \m I
        \end{pmatrix}\succeq \m 0.$$
        Finally, by applying the congruent transformation of matrix $\m \Theta$, we obtain equation \eqref{eq:th_LMI_input_feasibility}.
\end{proof}
}{}




The extension presented in \Cref{th:feasibility_linearization} builds upon \cite[Theorem~2]{egidio2022state} by incorporating the linearization error while maintaining convexity when considering the parameters of the initial ellipsoids $\set{\xi}$ as optimization variables. It is worth noting that in \cite[Theorem~2]{egidio2022state}, the problem is no longer an LMI if we consider $\m P$ as a variable. However, as demonstrated by \Cref{th:feasibility_linearization}, the optimization problem can be rendered convex through a congruent transformation when considering the square root of $\m P^{-1}$ as the variable of interest.
If $f$ is an affine function, the LMIs from \Cref{th:feasibility_linearization} can be simplified to the ones presented in~\cite[Theorem~2]{egidio2022state} by applying a congruent transformation.

\begin{figure}
    \centering
    \includegraphics[width=\linewidth]
    {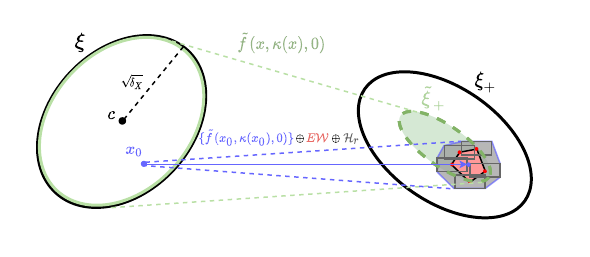}
    \caption{Local transition from ellipsoid $\set{\xi}$ to ellipsoid $\set{\xi}_+$.
    The affine controller $\kappa$ enforces that $\forall \vect x_0\in \set{\xi}$ and $\forall \vect w\in \set{W}$:
    $f(\vect x_0,\kappa(\vect x_0),\vect w)\in \{\tilde{f}(\vect x_0,\kappa(\vect x_0),\vect 0)\}\oplus \m E\set{W}\oplus \set{H}_r\subseteq \set{\xi}_+$
      where $\set{W}$ is the polytopic noise and $\set{H}_r=\hyperrectangle{\vect 0}{\tfrac{1}{2}\vect L r^2}$ is the noise resulting from the linearization \eqref{eq:inclusion_linearized_conservative} where $r^2=\delta_X+\delta_U+\delta_W$.
      The set $\tilde{\set{\xi}}_+ = q(\set{\xi})$ with $q(\vect x) = \tilde{f}(\vect x,\kappa(\vect x),\vect 0)$ (see \eqref{eq:minkowski-sum}) is an ellipsoid since an affine transformation of an ellipsoid is also an ellipsoid and the set $\m E \set{W}$ is a polytope since the linear transformation of a polytope is a polytope.}
    \label{fig:local_transition}
\end{figure}

We can visualize the terms of the Minkowski sum in~\eqref{eq:minkowski-sum} for our specific  setting in~\Cref{fig:local_transition}.
\keep{
 As the diameter of $\set{\xi}$ increases (i.e., $2\sqrt{\delta_X}$), the half-lengths of $\set{H}_r$ also increase. Consequently, the controller $\kappa$ must become more aggressive to ensure that the diameter of $\tilde{\set{\xi}}_+ \defEqual q(\set{\xi})$, defined after~\eqref{eq:minkowski-sum}, decreases. However, this might cause an increase in the magnitude of the control input, yielding a larger value of $\delta_U$ (see~\Cref{fig:local_transition_input}), which, in turn, amplifies the half-lengths of~$\set{H}_r$.
}
{
}

\subsection{Optimal control}
Since our aim is not only to design $\kappa$ to minimize a cost function, but also to optimize the shape $\m P$ of the starting ellipsoid $\set{\xi}$, we introduce a performance objective that involves minimizing a cost function \eqref{eq:cost_function} and maximizing the hypervolume of the preceding ellipsoid $\set{\xi}$.
As a reminder~\cite[Section 2.2.4]{boyd1994linear}, we can maximize the volume of an ellipsoid $\set{\xi}=\ellipsoid{\vect c}{\m P}$ by minimizing the convex function $-\log(\det(\m L))$ where $\m L$ is the square root of $\m P^{-1}$, i.e., $\m P^{-1}=\m L^{2}$.

The correctness and efficiency of the global RRT* algorithm~(\Cref{algo:RRTAbstraction}) are guaranteed by the following result.
\begin{corollary}\label{th:optimality_linearization}
\revise{For a given point $c$ and target set $\set{\xi}_+=\ellipsoid{\vect c_+}{\m P_+}$}, with a stage cost function $\set{J}$~\eqref{eq:cost_function} where $\m Q = \m S^\top \m S$ and $\lambda\in [0,1]$, the solution of the convex optimization problem
\begin{equation}\label{eq:biobjective}
\inf\limits_{\substack{\m L\succ \m 0,\m F, \revise{l}, \delta_X\ge 0,\delta_U\ge 0,\\\phi\ge 0,\beta_{ij}\ge 0,\tau_k\ge 0,\gamma\ge 0,\widetilde{\set{J}}}}\ \lambda \widetilde{\set{J}} +(1-\lambda) (-\log(\det(\m L)))
\end{equation}
\begin{align}
   \revise{\rm{s.t.}} \ & \eqref{eq:th_LMI_state_bound}, \eqref{eq:th_LMI_input_bound}, \eqref{eq:th_LMI_transition}, \eqref{eq:th_LMI_input_feasibility},\nonumber\\
    &\begin{pmatrix}
    \gamma \m I & \vect 0 & [\m L, \m F^\top, \vect 0]\m S^\top \\
    \DotSym & \widetilde{\set{J}}-\gamma & [\vect c^\top, \vect l^\top, \vect 1]\m S^\top  \\
    \DotSym & \DotSym & \m I
    \end{pmatrix}\succeq \m 0
    \label{eq:eq5:prob2}
\end{align}
satisfies
\begin{equation}\label{eq:cost_bound}
\widetilde{\set{J}} \ge \max_{\vect x\in \set{\xi}} \ \set{J}(\vect x,\kappa(\vect x))
\end{equation}
\revise{where the controller $\kappa(\vect x)\defEqual\m K(x-c) +l\in\set{U}$ with $\m K = \m F \m L^{-1}$ ensures a transition from $\set{\xi}=\ellipsoid{c}{\m P}$ with $\m P=\m L^{-2}$ to $\set{\xi}_+$, as in~\Cref{th:feasibility_linearization}}.
\end{corollary}
\keep{
\begin{proof} Same as in \cite[Corollary~1]{egidio2022state}.
\end{proof}
}{}
\keep{Note that unlike~\cite[Corollary 1]{egidio2022state}, the inequality~\eqref{eq:cost_bound} is generally not tight.}{}

\keep{
The parameter $\lambda$ governs the weight assigned to each criterion. Increasing its value will prioritize cost minimization, leading to a preference for reducing the volume of~$\set{\xi}$. Conversely, decreasing $\lambda$ will prioritize maximizing the volume of the starting ellipsoid~$\set{\xi}$, requiring larger control gains while mapping all states of~$\set{\xi}$ into~$\set{\xi}_+$. 
This also represents an exploitation/exploration trade-off, as maximizing the volume allows us to explore a larger portion of the state space, and minimizing the cost provides better (finer) solutions for already explored areas.
Because of that, $\lambda$ can be chosen adaptively throughout the execution.
}{
The parameter $\lambda$ governs the weight assigned to each criterion and also represents an exploitation/exploration trade-off, as maximizing the volume allows us to explore a larger portion of the state space, while minimizing the cost provides better (finer) solutions for already explored areas.
Because of that, $\lambda$ can be chosen adaptively throughout the execution.
}
\revise{A comprehensive study on how to fine-tune this meta-parameter is beyond the scope of this paper. However, a basic approach is to start with a smaller $\lambda$ to encourage designing larger cells (potentially leading faster to a feasible controller), and then increase its
value to favor smaller cells (with eventually lower transition costs).}

\keep{
\begin{figure}
    \centering
    \includegraphics[width=0.6\linewidth] {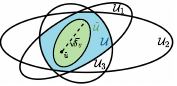}
    \caption{
    Input space of the system $\Sys$ (blue) and the set of inputs actually used by the controller $\kappa$ on the ellipsoid $\xi$ (green). The input space is denoted by $\set{U} = \set{U}_1\cap \set{U}_2\cap \set{U}_3$ where $\set{U}_1,\set{U}_2,\set{U}_3$ are ellipsoids, and $\widetilde{\set{U}} = \kappa(\set{\xi})$ represents the set of inputs used by $\kappa$ on $\set{\xi}$. 
    The controller $\kappa$ is designed such that $\widetilde{\set{U}}\subseteq \set{U}$.
    }
    \label{fig:local_transition_input}
\end{figure}
}{}

\section{Numerical experiments}\label{sec:experiments}
\begin{figure}
    \centering
    \includegraphics[width=\linewidth]
    {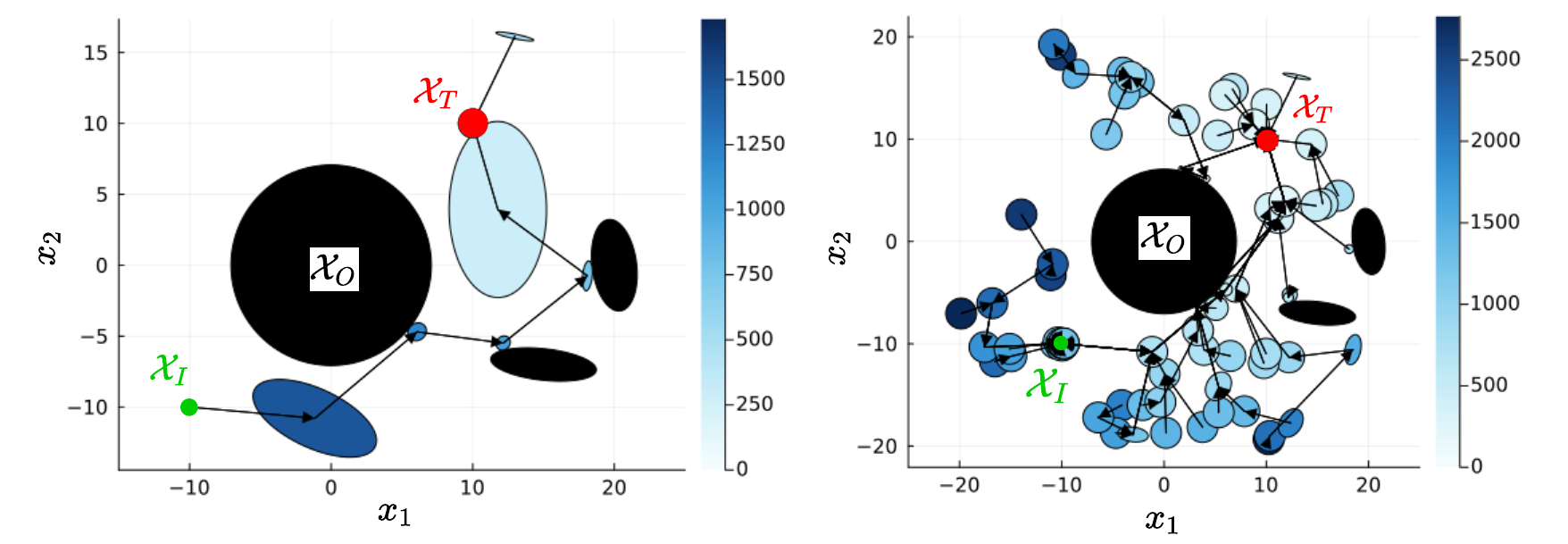}    
    \caption{    
    (\Cref{sec:experiments}) 
    Solution obtained from~\Cref{algo:RRTAbstraction} (left) and from the improved RRT* version (right).
     The initial set $\set{X}_I$, target set $\set{X}_T$ and obstacles $\set{X}_O$ are respectively in \revise{green, red and black}.
    The color map illustrates the value function $v$ for $\Sys$.
    Left: Construction of $\widetilde{\Sys}_1$ and $v_1$ stops as soon as a feasible solution is found. Right: Continue to extend and improve the abstraction $\widetilde{\Sys}_2$ and the value function $v_2$. We added a bound on the cell volume to easily visualize the entire tree.}
    \label{fig:global_algorithm}
\end{figure}

\keepSecond{
\subsection{One single transition}\label{sec:singleTransition-planar}
In this example we study aspects of determining a single transition for a given target set $\set{\xi}_+ = 
\ellipsoid{\vect c_+}{\m P_+}$ and initial point $\vect c$ with 
$$\vect c_+ = \begin{pmatrix}
4\\ 4
\end{pmatrix}, 
\ \m P_+ = 
\begin{pmatrix}
2 & 0.2\\
0.2 & 0.55
\end{pmatrix}, \ \vect c = \begin{pmatrix}
1\\
1
\end{pmatrix}.$$
Consider the nonlinear system \eqref{eq:bounded_disturbances} given by 
\begin{equation}\label{eq:2d-system}
f(\vect x, \vect u, \vect w) = 
\begin{pmatrix}
    1.1 x_1 - 0.2 x_2 - \rho  x_2^3 + u_1 + w_1\\
    0.2 x_1 + 1.1 x_2 + \rho  x_1^3 + u_2 + w_2
\end{pmatrix}
\end{equation}
where $\rho\ge 0$. The control input $\vect u$ is constrained by the set $\set{U} = \set{U}_1 \cap \set{U}_2 \cap \set{U}_3$ with 
$$
\footnotesize{
\set{U}_1 = 
\operatorname{H}\left(\vect 0, \begin{pmatrix}
4\\5
\end{pmatrix}
\right), \set{U}_2 = \ball{\vect 0}{5}, \set{U}_3 = \operatorname{H}\left(\vect 0,
\begin{pmatrix}
0.05 & 0\\
0 & 0.033
\end{pmatrix}\right)}
$$
and the exogenous input $\vect w$ by the set $\set{W} = \hyperrectangle{\vect 0}{
(\omega_{\text{max}}, \ \omega_{\text{max}})^\top
}$ with $\omega_{\text{max}}\ge 0$. 
The level of non-linearity in the system is determined by the parameter $\rho$, where $\rho=0$ corresponds to an affine dynamical system. On the other hand, the presence of noise is controlled by the parameter $\omega_{\text{max}}$, where $\omega_{\text{max}}=0$ corresponds to a deterministic system.

For several different values of $\rho$, $\omega_{\text{max}}$ and $\lambda$, we solved the optimization problem of~\Cref{th:optimality_linearization} considering the quadratic cost function~\eqref{eq:cost_function} with $\m Q = \text{diag}(\m I, \ \m I, \ 1)$, i.e., 
$\set{J}(\vect x,\vect u) = \vect x^\top \vect x + \vect u^\top \vect u + 1$.

On average, each solution to \Cref{th:optimality_linearization} was found in $0.0176$ seconds on an Intel\textsuperscript{\textregistered} Core\texttrademark{} i7-10610U CPU $1.80$ GHz $\times 8$ with $16$ GB of memory and using the Julia JuMP \cite{dunning2017jump} interface with the Mosek solver on Windows $10$.

Firstly, observing the left-hand side of~\Cref{fig:single-transition:1}, we can deduce that while $\widetilde{\set{J}}$ is an upper bound for the cost of the transition from $\set{\xi}$ to $\set{\xi}_+$, in practice this cost is considerably lower for a significant part of the cell.
As expected, by comparing~\Cref{fig:single-transition:1} and~\Cref{fig:single-transition:2}, the worst case cost $\widetilde{\set{J}}$ of the controller decreases as $\lambda$ increases. This is achieved by reducing both the size of the initial ellipsoid $\xi$ and the maximum magnitude of the inputs actually used by the controller $\kappa$ on $\xi$.
 Finally, comparing the results of \Cref{fig:single-transition:1} and \Cref{fig:single-transition:3}, we notice that when we increase the non-linearity ($\rho$) and the noise bound ($\omega_{\text{max}}$), the controller has to produce a smaller ellipsoid $\tilde{\set{\xi}}_+$ more centered in $\set{\xi}_+$, and to do so, it must reduce the size of the initial ellipsoid.
 
\begin{figure}
    \centering
    \includegraphics[width=0.49\textwidth]{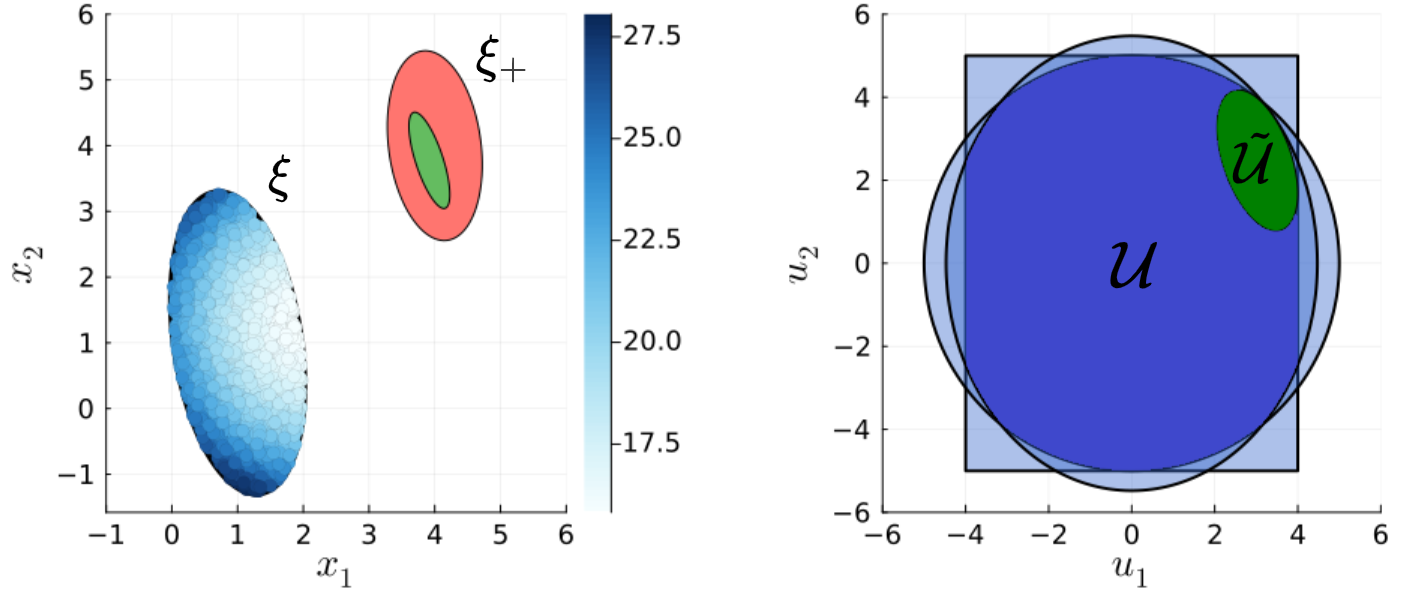}
    \caption{(\Cref{sec:singleTransition-planar})
        Solution provided by~\Cref{th:optimality_linearization} for $\omega_{\text{max}} = 0.1$, $\rho = 0.0005$ and $\lambda = 0.01$. 
          Left: Value of the cost function $\set{J}(x,\kappa(x))$ for the closed-loop system (color map) and the ellipsoids $\set{\xi}_+$ (red) and $\tilde{\set{\xi}}_+$ (green). 
          Right: The sets $\set{U}_1$, $\set{U}_2$, and $\set{U}_3$ (blue), and their intersection $\set{U}$ (dark blue). 
          The set $\widetilde{\set{U}} \defEqual \kappa(\set{\xi})$ (green) represents the inputs actually employed by the controller $\kappa$ on $\xi$. 
          The upper bound on the transition cost $\widetilde{\set{J}} = 28.1$ and the volume of the starting ellipsoid $\vol(\set{\xi}) = 6.65$.
    }
    \label{fig:single-transition:1}
\end{figure}

\begin{figure}
    \centering
    \includegraphics[width=0.49\textwidth]{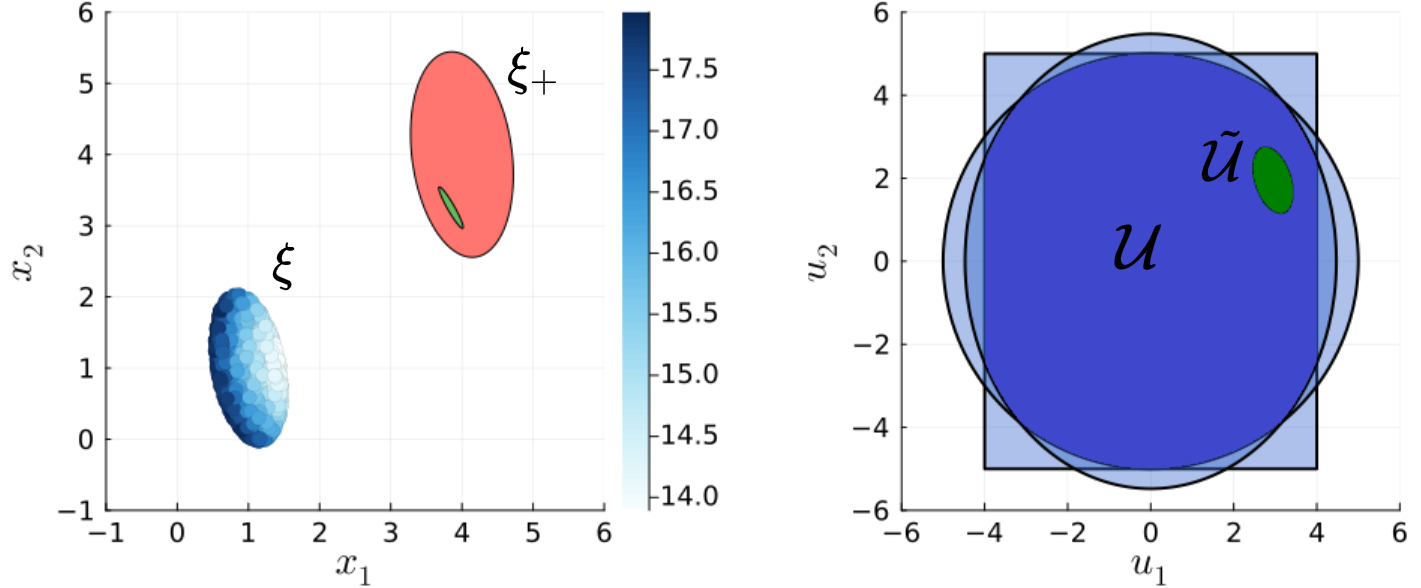}
    \caption{
        (\Cref{sec:singleTransition-planar}) Solution provided by~\Cref{th:optimality_linearization} for $\omega_{\text{max}} = 0.1$, $\rho = 0.0005$ and $\lambda = 0.3$.   We have $\widetilde{\set{J}} = 18$, $\vol(\set{\xi}) =  1.44$.
    }
    \label{fig:single-transition:2}
\end{figure}

\begin{figure}
    \centering
    \includegraphics[width=0.49\textwidth]{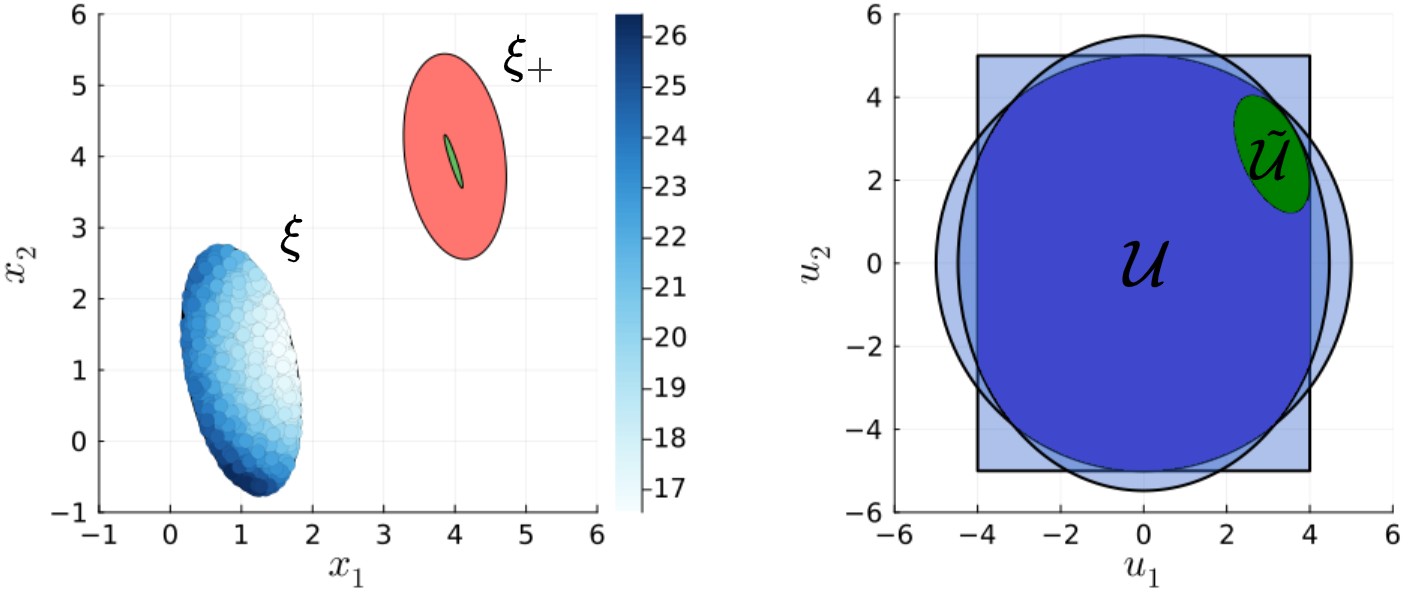}
    \caption{
        (\Cref{sec:singleTransition-planar}) Solution provided by~\Cref{th:optimality_linearization} for $\omega_{\text{max}} = 0.15$, $\rho = 0.001$ and $\lambda = 0.01$. 
          We have $\widetilde{\set{J}} = 26.56$, $\vol(\set{\xi}) = 3.82$.
    }
    \label{fig:single-transition:3}
\end{figure}

\begin{figure}
    \centering
    \includegraphics[width=0.5\textwidth]{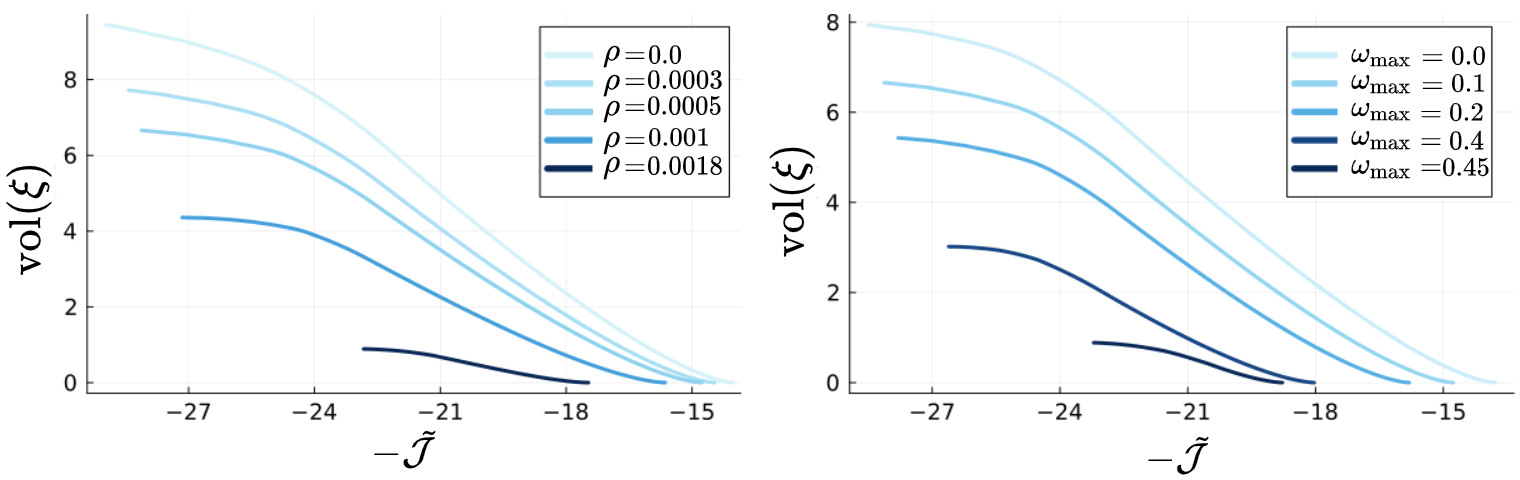}
    \caption{
       (\Cref{sec:singleTransition-planar}) Pareto front of the bi-objective function \eqref{eq:biobjective} for different values of $\rho$ and $\omega_{\text{max}}$. 
     Each curve decreases as the value of $\lambda$ increases from $0$ to~$1$. Left: $\omega_{\text{max}} = 0.1$. Right: $\rho = 0.0005$.
    }
    \label{fig:single-transition:part2-1}
\end{figure}

\Cref{fig:single-transition:part2-1} illustrates that as the values of $\omega_{\text{max}}$ and $\rho$ increase, the Pareto front of the bi-objective function~\eqref{eq:biobjective} consistently exhibits worse outcomes, i.e.,  higher costs and smaller volumes of the initial ellipsoid, across all values of the weighting parameter $\lambda$. 
\subsection{Optimal control}\label{sec:optimal_control_problem}
}{}
\keep{In this example we provide one possible application of the approach for the optimal control of the two-dimensional dynamical system $\set{S}$ introduced in~\Cref{sec:singleTransition-planar} with $\rho=0.005$, $\omega_{\max} = 0.1$, and the same stage cost function, which penalizes states and inputs that are distant from the origin.}
{
    We consider the nonlinear dynamical system \eqref{eq:bounded_disturbances} given by 
    \begin{equation}\label{eq:2d-system}
    f(\vect x, \vect u, \vect w) = 
    \begin{pmatrix}
        1.1 x_1 - 0.2 x_2 - \rho  x_2^3 + u_1 + w_1\\
        0.2 x_1 + 1.1 x_2 + \rho  x_1^3 + u_2 + w_2
    \end{pmatrix}
    \end{equation}
    where $\rho=0.0005$. The control input $\vect u$ is constrained by the set $\set{U} = \set{U}_1 \cap \set{U}_2 \cap \set{U}_3$ with 
    $$
    \footnotesize{
    \set{U}_1 = 
    \operatorname{H}\left(\vect 0, \begin{pmatrix}
    4\\5
    \end{pmatrix}
    \right), \set{U}_2 = \ball{\vect 0}{5}, \set{U}_3 = \operatorname{H}\left(\vect 0,
    \begin{pmatrix}
    0.05 & 0\\
    0 & 0.033
    \end{pmatrix}\right)}
    $$
    and the exogenous input $\vect w$ by the set $\set{W} = \hyperrectangle{\vect 0}{
    (\omega_{\text{max}}, \ \omega_{\text{max}})^\top
    }$ with $\omega_{\text{max}}=0.01$.
    The specification $\Sigma$ is clear from~\Cref{fig:global_algorithm} and we consider the quadratic cost stage cost function~\eqref{eq:cost_function} with $\m Q = \text{diag}(\m I, \ \m I, \ 1)$, i.e., $\set{J}(\vect x,\vect u) = \vect x^\top \vect x + \vect u^\top \vect u + 1$.
}

A first feasible solution $\widetilde{\Sys}_1$ (see~\Cref{fig:global_algorithm}) was found after $18$ seconds with only $8$ cells and $7$ transitions created\keep{, in the same computational setup as in the previous example.}{.} We can continue to expand the tree to explore other paths in the state-space as illustrated with abstraction $\widetilde{\Sys}_2$ in~\Cref{fig:global_algorithm} (right).
As anticipated, the controller exhibits a preference for solutions that pass near the origin while circumventing the obstacle.

Given $x_0=(-10,-10)\in \set{X}_I$, the guaranteed total cost by $\widetilde{\Sys}_1$ is $v_1(x_0) = 1732$ whereas the true total cost
of this specific trajectory is $1337$.
For $\widetilde{\Sys}_2$, the guaranteed total cost is $v_2(x_0) = 992$ 
whereas the true total cost is $875$.
The control cost guaranteed by $\widetilde{\Sys}_2$ is better than by $\widetilde{\Sys}_1$ for two main reasons: 1) we continued to improve the abstraction using the RRT* variant, and 2) because by considering smaller cells, the cost of transitioning from one cell to another in the worst case is lower. Nevertheless, a similar part of the state space is covered by both abstractions, which is achieved with far fewer cells for $\widetilde{\Sys}_1$. 
%


\section{Conclusion}
We provided a tractable algorithm for the optimal control of $L$-smooth nonlinear dynamical systems.
Firstly, the proposed algorithm circumvents the necessity of discretizing the input space by employing a set of local feedback controllers. This is done to ensure deterministic transitions, thereby eliminating the non-determinism associated with abstraction, a common limitation in classic approaches.
Secondly, the combined use of ellipsoid-based covering and affine local controllers leverages the power of LMIs and convex optimization,
allowing the creation of larger and non-standard cells.
Thirdly, the use of a lazy approach and non-uniform cells significantly reduces the complexity of the abstraction, i.e., the number of abstract states, when addressing a specific control problem.

As future work, we plan to demonstrate the efficiency of this approach on higher-dimensional nonlinear dynamical systems.

\bibliographystyle{IEEEtran}
\bibliography{ref}

\end{document}